\documentclass{birkmult}
\usepackage{graphics}
\usepackage{amssymb}
\def\Z{{\mathbb Z}}
\def\sym{{\mathcal S}}

 \newtheorem{thm}{Theorem}[section]
 \newtheorem{cor}[thm]{Corollary}

 \theoremstyle{definition}
 \newtheorem{defn}[thm]{Definition}
 \theoremstyle{remark}

 \numberwithin{equation}{section}
\begin{document}
\title{Automorphism groups of circulant graphs --- a survey}
\author{Joy Morris}
\address{Department of Mathematics and Computer Science \\
University of Lethbridge \\
Lethbridge, AB. T1K 6R4. Canada}
\email{joy@cs.uleth.ca}
\thanks{The author gratefully acknowledges support from the National Science and Engineering Research Council of Canada (NSERC).}
\begin{abstract}
A circulant (di)graph is a (di)graph on $n$ vertices that admits a cyclic automorphism of order $n$. This paper provides a survey of the work that has been done on finding the automorphism groups of circulant (di)graphs, including the generalisation in which the edges of the (di)graph have been assigned colours that are invariant under the aforementioned cyclic automorphism.
\end{abstract}
\subjclass{05C25}
\keywords{circulant graphs, automorphism groups, algorithms}
\maketitle
\section{Introduction}

The aim of this paper is to provide a history and overview of work that has been done on finding the automorphism groups of circulant graphs.  We will focus on structural theorems about these automorphism groups, and on efficient algorithms based on these theorems, that can be used to  determine the automorphism group of certain classes of circulant graphs.

We must begin this discussion by defining the terms that will be central to the topic.  In what follows, we sometimes refer simply to ``graphs" or to ``digraphs," but all of our definitions, and many of the results on automorphism groups, can be generalised to the case of  ``colour digraphs:" that is, digraphs whose edges not only have directions, but colours.
\begin{defn}
Two graphs $X=X(V,E)$ and $Y=Y(V',E')$ are said to be 
{\bf isomorphic} if there is a bijective mapping $\phi$ from the vertex set $V$ to the 
vertex set $V'$ such that $(u,v)\in E$ if and only if $(\phi(u),\phi(v))
\in E'$.  The mapping $\phi$ is called an {\bf isomorphism}.
We denote the fact that $X$ and $Y$ are isomorphic by $X
\cong Y$.
\end{defn}
That is, an isomorphism between two graphs is a bijection on the vertices that preserves edges and nonedges.  In the case of digraphs, an isomorphism must also preserve the directions assigned to the edges, and in the case of colour (di)graphs, the colours must also be preserved.

This definition has the following special case:
\begin{defn}
An {\bf automorphism} of a (colour) (di)graph is an isomorphism from the (colour) (di)graph to itself.
\end{defn}

When we put all of the automorphisms of a (colour) (di)graph together, the result is a group:
\begin{defn}
The set of all automorphisms of a graph $X$ forms a group, denoted $\rm{Aut}(X)$, the {\bf automorphism group} of $X$.
\end{defn}

Now that we know what automorphism groups of graphs are, we must define circulant graphs.
\begin{defn}
A {\bf circulant graph} $X(n;S)$ is a Cayley graph on $\Z_n$.  That is, it is a graph whose vertices are labelled $\{0,1,\ldots, n-1\}$, with two vertices
labelled $i$ and $j$ adjacent iff $i-j \pmod{n} \in S$, where $S \subset \Z_n$
has $S=-S$ and $0 \not\in S$.
\end{defn}
For a circulant digraph, the condition that $S=-S$ is removed.  For a colour circulant (di)graph, each element of $S$ also has an associated (not necessarily distinct) colour, which is assigned to every edge
whose existence is a consequence of that element of $S$.

With these basic definitions in hand, we can explore the history of the question: what do the automorphism groups of circulant graphs look like, and how can we find them?

Although we define most of the terms used in this paper, for permutation group theoretic terms that are not defined, the reader is referred to Wielandt's book on permutation groups \cite{W}, which has recently come back into print in his collected works \cite{W3}.  Another good source is Dixon and Mortimer's book \cite{DM}.  Throughout this paper, the symmetric group on $n$ points is denoted by $\sym_n$.

\section{History}
In 1936, 
K\"{o}nig \cite{Konig} asked the following question:
``When can a given abstract group be interpreted as the group
of a graph and if this is the case, how can the corresponding graph
be constructed?  This same question could be asked for directed
graphs." (This quote is from the English translation published by Birkha\"{u}ser in 1990).  By the group of a graph, he is referring to the automorphism group.

This question was answered by
Frucht, in 1938 \cite{Frucht}.  The answer was yes; in fact, it went further: there are infinitely many such graphs for any group $G$.

By no means was this the end of the matter.  One method of construction involved the creation of graphs on large numbers of vertices, that encoded the colour information from the Cayley colour digraph of a group, into structured subgraphs.  This construction led to graphs that had, in general, many more vertices than the order of the group.

One major area of research that spun off from this, was the search for ``graphical regular representations" of a particular group $G$: that is, graphs whose automorphism group is isomorphic to $G$, and whose number of vertices is equal to the order of $G$.  This is not a topic that we will pursue further in this paper, however.

A related question was presumably also considered by mathematicians.  That is,
given a particular representation of permutation group $G$, is there a graph $X$ for which Aut($X$) $\cong G$ as permutation groups? We know of no reference for this question prior to 1974, however \cite{Brian2}.

The answer to this question did not prove to be so straightforward.
For example, if $\rho = (0 \ 1 \ 2 \ldots n)$, no graph has $\langle \rho \rangle$
as its automorphism group, because the reflection that maps $a$ to $-a$ for every $0\le a\le n$ will always be an automorphism of any such graph.
Another way of explaining this, is that the dihedral group acting on $n$ elements has the same orbits on unordered pairs $(i,j)$ as $Z_n$ has.  This shows that there is not always a graph $X$ whose automorphism group is a particular representation of a permutation group.

Although it is not our aim to discuss this version of the question, we will give here the most significant results that have been obtained on it.

The first result is due to Hemminger \cite{Hemminger}, in 1967.
\begin{thm}
Let $G$ be a transitive abelian permutation group that is abstractly isomorphic to $\Z_{n_1}\times \ldots \times \Z_{n_k}$.  Then if the number of factors in the direct product of order 2 is not 2, 3, or 4, and the number of factors of order 3 is not 2, there is a directed graph whose automorphism group is isomorphic to $G$.
\end{thm}

In 1981, Godsil \cite{Godsil2} proved the following result.
\begin{thm}
Let $G$ be a finite permutation group.  A necessary condition for the Cayley graph $X=X(G;S)$ to have $G$ as its automorphism group, is that the subgroup of $G$ that fixes some vertex of $X$, is isomorphic to the automorphisms of $G$ that fix $S$ set-wise.  
\end{thm}
He also proved that this condition is sufficient for many $p$-groups, and obtained necessary and sufficient conditions for the occurrence of the dihedral groups of order $2^k$, and of certain Frobenius groups, as the full automorphism groups of vertex-transitive graphs and digraphs. 

In 1989, a paper by Zelikovskij \cite{Zeli} appeared in Russian.  We are unaware of a translation, so can only report the main result as stated in the English summary: that for every finite abelian permutation group $G$ whose order is relatively prime to 30, the paper provides necessary and sufficient conditions for the existence of a simple graph whose automorphism group is isomorphic to $G$.  Note that the facts that Zelikovskij produces a simple graph and apparently does not require transitivity, represent improvements over Hemminger's result, although he does not cover all of the orders that Hemminger does.

In 1999, Peisert proved \cite{Peisert} the following result.
\begin{thm}
If two permutation groups each have representations as automorphism groups of graphs, then the direct product of these representations, will also have a graph for which it is isomorphic to the automorphism group, unless the two original groups are isomorphic as permutation groups, transitive, and have a unique graph for which they are the automorphism group (up to isomorphism).
\end{thm}

The approach this question that we will follow in this paper, is to ask what is the automorphism group of a given graph $X$?  Again, the answer is not easy, so we will limit our consideration to circulant graphs.

\section{Algorithms for finding automorphism groups}

As one focus of this paper will be finding efficient algorithms for calculating the automorphism group of a graph, it is necessary to spend some time considering what makes an algorithm ``efficient" in this regard.

There is, after all, a very straightforward algorithm that is guaranteed to find the automorphism group of any graph on $n$ vertices: simply consider every possible permutation in $\sym_n$, and include those that turn out to be automorphisms of our graph.  As $\sym_n$ has order $n!$, this algorithm is exponential in $n$.

Depending on our goal, it may not be possible to do any better than this.  Specifically, if our goal is to list every automorphism in the automorphism group, then as there may be as many as $n!$ elements, generating the list in anything less than exponential time will not be feasible in general. 

If the number of prime factors of $n$ is bounded  with a sufficiently low bound, then we actually can do more than this.  In such cases, the number of subgroups of $\sym_n$ that can be automorphism groups of circulant graphs may be sufficiently small that we can choose the automorphism group of our graph from a short list.  We will not be listing all of its elements, but we may be able to describe its structure precisely without the use of generating sets, as for example by saying it is isomorphic to $\sym_p \times \sym q$.

However, when the number of prime factors is unbounded, we lose this ability to explicitly describe the automorphism group.  Thus, when we consider algorithms intended to determine the automorphism group of a circulant graph on a number of vertices whose number of factors is unbounded, what we are looking for is not an explicit listing of all of the elements of the group, but the provision of a generating set for the group.  

Each of the algorithms that we provide below, runs in polynomial time in the number of vertices of the graph, and provides either a generating set for the automorphism group of the graph, or (where possible) a precise structural description of the automorphism group.

\section{Circulant graphs on a prime number of vertices}

In 1973, Alspach proved the following result on the automorphism groups of circulant graphs on a prime number of vertices.

\begin{thm}\cite{Brian}
Let $p$ be prime.
If $S= \emptyset$ or $S = \Z_p^*$, then Aut($X$)$=\sym_p$; otherwise,
Aut($X$)$=\{T_{a,b}: a \in E(S), b \in \Z_p\}$, where $T_{a,b}(v_i)=v_{ai+b}$,
and $E(S)$ is the largest even-order subgroup of $\Z_p^*$ such that $S$ is
a union of cosets of $E(S)$.
\end{thm}
Notice that since $S=-S$, $S$ must be a union of cosets of $\{1,-1\}$, so $E(S)$ can always be found.

This leads to the following algorithm for finding the automorphism group of such a graph.
\smallskip

\noindent
{\bf Algorithm for finding Aut($X$)}:
\begin{enumerate}
\item If $S=\emptyset$ or $S=\Z_p^*$, Aut($X$)$=\sym_p$.
\item For each even-order subgroup $H$ of $\Z_p^*$, verify whether $S$ is a
union of cosets of $H$.
\item Since $\Z_p^*$ is cyclic, there is one $H$ of every order dividing $p-1$.  Set $E(S)$ to be the largest $H$ that satisfies (2).
\item If Aut($X$)$\neq \sym_p$, Aut($X$)$=\{T_{a,b}:a \in E(S), b\in \Z_p\}$.
\end{enumerate}

Proof of Alspach's theorem relies on the following theorem by Burnside.

\begin{thm}\cite{burn1}
If $G$ is a transitive group acting on a prime number $p$ of elements, then 
either $G$ is doubly transitive or $G=\{T_{a,b}: a \in H <\Z_p^*, b\in \Z_p\}$.
\end{thm}

Although the algorithm given above may be the most natural way to create an algorithm from the statement of Alspach's theorem, we will re-state the algorithm slightly differently.  The purpose of this, is to provide a closer parallel to the algorithms that we will be constructing subsequently, to cover other possible numbers of vertices.
\smallskip

\noindent
{\bf Alternate algorithm}:
\begin{enumerate}
\item Find $A$, the set of all multipliers $a \in \Z_p^*$ for which $aS=S$.
\item If $A=\Z_p^*$, then Aut($X$)$=\sym_p$.
\item Otherwise, Aut($X$)$=\{T_{a,b}: a\in A, b\in \Z_p\}$.
\end{enumerate}

Before presenting any generalisations of this algorithm, the concept of wreath products will be required.

\section{Wreath products}

Although we will present a formal definition of the wreath product of two graphs in a moment, we will first give a description which may make the formal presentation easier to follow.  If we are taking a wreath product of two (di)graphs, $X$ and $Y$, we replace every vertex of $X$ by a copy of the (di)graph $Y$.  Between two copies of $Y$, we include all edges (or all arcs in a particular direction) if there was an edge (or an arc in the appropriate direction) between the corresponding vertices of $X$.  Now for the formal definition.

\begin{defn}
The wreath product of the (di)graph $X$ with the (di)graph $Y$, denoted $X \wr Y
$,
is defined in the following way.

The vertices of $X \wr Y$ are the ordered pairs $(x,y)$ where $x$ is a vertex of
$X$ and $y$ is a vertex of $Y$.  There is an arc (or edge) from the vertex $(x_1
,y_1)$
to the vertex $(x_2,y_2)$ if and only if one of the following holds:
\begin{enumerate}
\item
$x_1=x_2$ and $(y_1,y_2)$ is an arc (or edge) of $Y$; or
\item
$(x_1,x_2)$ is an arc (or edge) of $X$.
\end{enumerate}
\end{defn}

We say two sets of vertices $A$ and $B$ are {\bf wreathed} if either
$ab$ is an edge for every $a\in A$ and every $b\in B$, or $ab$ is a nonedge
for every $a\in A$ and every $b\in B$.

In graph theory, what we have called the ``wreath" product is also often called the ``lexicographic" product, and has also been called the ``composition" of graphs.
The notion of a wreath product is also defined on groups; in fact, the term ``wreath product" comes from group theory.  Although it can be defined on abstract groups, we will only be considering permutation groups, where the definition is simpler, so it is this definition that we provide.

\begin{defn}
The wreath product of two permutation groups, $H$ and $K$, acting on sets
$U$ and $V$ respectively, is 
the group of all permutations $f$ of $U \times V$ for
which there exist $h \in H$ and an element $k_u$ of $K$ for each $u \in U$
such that $$f((u,v))=(h(u),k_{h(u)}(v))$$ for all $(u,v)\in U \times V$.
It is written $H \wr K$.
\end{defn}

It is easy to verify that 
$$\rm{Aut}(X)\wr\rm{Aut}(Y)\le\rm{Aut}(X\wr Y)$$
is true for any graphs $X$ and $Y$; in fact, it is often the case that equality holds.
 
\section{Circulant graphs on $pq$ vertices, or $p^n$ vertices}

The following theorem, proven by Klin and P\"{o}schel in 1978, characterises the automorphism groups of circulant graphs on $pq$ vertices, where $p$ and $q$ are distinct primes.

\begin{thm}\cite{KPa}
If $G$ is the automorphism group of a circulant graph on $pq$ vertices, then $G$ is one of:
\begin{enumerate}
\item $\sym_{pq}$;
\item $A_1 \wr A_2$, or $A_2 \wr A_1$, where $A_1$ and $A_2$ are automorphism
groups of circulant graphs on $p$ and $q$ vertices, respectively; 
\item $\sym_p \times A_2$, or $A_1\times \sym_q$, where $A_1$ and $A_2$ are automorphism groups of
circulant graphs on $p$ and $q$ vertices, respectively; or,
\item $\{T_{a,b}: a \in A \le \Z_n^*, b \in Z_n\}$ (a subgroup of the holomorph).
\end{enumerate}
\end{thm}

The following algorithm can be constructed from this characterisation, to determine the automorphism group of a graph on $pq$ vertices.
\smallskip

\noindent
{\bf Algorithm}:
\begin{enumerate}
\item If $S=\emptyset$ or $S=\Z_{pq}-\{0\}$, then Aut($X$)$=\sym_{pq}$. END.
\item If $(v_0\ \ v_p \ \ v_{2p} \ldots\ \ v_{(q-1)p)}\in$Aut($X$), then 
Aut($X$)$=A_1 \wr A_2$, where $A_1$ is the automorphism group of the induced subgraph of $X$ on 
the vertices $\{v_0,v_q,\ldots,v_{(p-1)q}\}$ and $A_2$ is the automorphism group of the induced subgraph of $X$ on the vertices $\{v_0,v_p,\ldots,v_{(q-1)p}\}$. Use the previous algorithm (for finding the automorphism group of a circulant graph on a prime number of vertices) to find $A_1$ and $A_2$. END.
\item Repeat (2) with the roles of $p$ and $q$ reversed.
\item Let $A$ be the group of all multipliers $a$ in $\Z_{pq}^*$ for which $aS=S$. 
\item Define $E_p$ by
$E_p=\{T_{a,b}:a \in A, a \equiv 1\pmod{p}, b\in q\Z_p\}$; if
$E_p \cong$AGL$(1,p)$ then Aut$(X)=\sym_p \times A_2$, where $A_2$ is as in step (2).  Use the previous algorithm to find $A_2$. END.
\item Repeat (4) with the roles of $p$ and $q$ reversed, and $A_1$ taking the role of $A_2$. END.
\item Otherwise, Aut($X$)$=\{T_{a,b}:a\in A, b\in \Z_{pq}\}$.
\end{enumerate}

The full automorphism group of circulant graphs on $p^n$ vertices has been determined, first by Klin and P\"{o}schel, and later independently by Dobson; both results are unpublished.

The full result is technical, but this nice result gets much of the way:

\begin{thm}\cite{Ted, KP2}
A circulant graph on $p^n$ vertices is either a wreath product, or its automorphism group has a 
normal Sylow $p$-subgroup.
\end{thm}

\section{The square-free case}

There is no known characterisation for the automorphism groups of circulant graphs in the general square-free case that is as straightforward as the results we have described above, when the number of vertices is $p$ or $pq$ and $p, q$ are distinct primes.

However, Dobson and Morris \cite{bteam1} did prove the following structural theorem about the automorphism groups of circulant graphs of arbitrary square-free order.  The definition of a group being 2-closed is quite technical, but the important thing to know about 2-closed groups is that the automorphism group of a vertex-transitive graph or digraph is always 2-closed.  The terms used in the second point of this theorem are defined later in this paper, in the section ``A Strategy," for readers who are interested.
\begin{thm}
Let $mk$ be a square-free integer and $G\le \sym_{mk}$ be
$2$-closed and contain a regular cyclic subgroup, $\langle \rho\rangle$.
Then one of the following is true:
\begin{enumerate}
\item $G = G_1\cap G_2$, where $G_1 = \sym_r\wr H_1$ and $G_2 = H_2\wr \sym_k$,
where $H_1$ is a $2$-closed group of degree $mk/r$, $H_2$ is a
$2$-closed group of order $m$, and $r\vert m$; or
\item there exists a complete block system ${\mathcal B}$ of $G$
consisting of $m$ blocks of size $k$, and there exists $H\triangleleft G$
such that $H$ is transitive, $2$-closed, and $\langle \rho \rangle \le H =
H_1\times H_2$ (with the canonical action), where $H_1\le
\sym_{m}$ is $2$-closed and $H_2\le \sym_{k}$ is $2$-closed and
primitive.
\end{enumerate}
\end{thm}

Unlike the previous structural results by Alspach, and by Klin and P\"{o}schel, this clearly does not provide a clear, short list of groups from among which the automorphism group of a circulant graph on any square-free number of vertices must be chosen.  However, it is sufficient to allow Dobson and Morris, in a subsequent paper \cite{algo}, to construct the following algorithm that will determine the automorphism group of any such graph.
\smallskip

\noindent
{\bf General algorithm ($n$ square-free)}

\noindent{\bf Inputs:} The number $n$ of vertices of $X$, and the
connection set $S\subseteq \Z_n$.
\begin{enumerate}
\item Let $A$ be the group of all multipliers $a$ in $\Z_n^*$ for which $aS=S$.
\item For each prime divisor $p$ of $n$, define $E_p$ by
$E_p=\{T_{a,b}:a \in A, a \equiv 1\pmod{p}, b\in \frac{n}{p}\Z_p\}$; if
$E_p \cong$AGL$(1,p)$ then define $E_p=\sym_p$.
\item Let $p_1, \ldots, p_t$ be all primes such that $E_p=\sym_p$.
\begin{enumerate}
\item For each pair of distinct primes $p_i, p_j$, if the transposition that switches $n/p_i+kp_ip_j$ with $n/p_j+kp_ip_j$ for every $0\le k\le n/p_ip_j$ and fixes all other vertices is an automorphism of $X$, then $E_{p_ip_j}=\sym_{p_ip_j}$.
Otherwise, $E_{p_ip_j}=E_{p_i}E_{p_j}$.
\item Define a relation $R$ on $\{p_1, \ldots, p_t\}$ by $p_i R p_j$ iff there
exists a sequence of primes $p_{k_1}, \ldots, p_{k_s}$ where $p_i=p_{k_1}, p_j=p_{k_s}$ such that $E_{\rm{lcm}(p_{k_l},p_{k_{l+1}})}=\sym_{\rm{lcm}(p_{k_l},p_{k_{l+1}})}$.  This is an equivalence relation.  For each equivalence class $\mathcal E_i$, let $m_i=\Pi_{j\in \mathcal E_i p_j}$, then let $E_{m_i}=\sym_{m_i}$.
\item For any divisor $m$ of $n$, $$E_m=\Pi_{\mathcal E}\sym_{\rm{gcd}(m_i,m)}
\Pi_{p\vert \rm{gcd}(n/p_1\ldots p_t,m)}E_p.$$
\end{enumerate}
\item For each composite divisor $m$ of $n$, let $A_m=\{a \in A: a \equiv 1($mod $\frac{n}{m})\}$, and let $A_n=A$.
Define $E_m'=\langle E_m, A_m\rangle$.  
\item Let $G=E_n'$.  For each complete block system of $E_n'$, ${\mathcal B}$, consisting of $n/k$ blocks of size $k$, do:
\begin{enumerate}
\item For each complete block system of $E_n'$, ${\mathcal D}$, consisting
of $n/kk'$ blocks of size $kk'$, determine whether or not
$$\rho^{n/k}\vert_{D_0},$$
the mapping that acts as $\rho^{n/k}$ on the vertices of $D_0$, and fixes all other vertices,
is an automorphism of $X$.  If it is, let our new $G$ be the group generated by the old $G$ together with every $E_k'\vert_D$, where $D\in {\mathcal D}$; if not, leave $G$ unchanged.
\end{enumerate}
\item $G$ is the automorphism group of $X$.
\end{enumerate}

To avoid some technical details, we have oversimplified some of the notation in this algorithm; the astute reader may observe, for example, that if we define $E_p=\sym_p$, then $E_p$ is acting on $p$ vertices, and subsequent products involving $E_p$ ($E_m=E_pE_q$, etc.,) may not be properly defined.  These issues are properly addressed in the paper that presents the algorithm, but as the intent of this paper was only to give the flavour of the algorithm, the technical details seemed likely to unnecessarily complicate our presentation.

In the same paper, Dobson and Morris prove that this algorithm runs in polynomial time on the number of vertices of the graph.

\section{A Strategy}

In this section, we outline a general strategy that is often used to obtain structural results about circulant graphs.

\begin{defn}
Let $V$ be a set, and $G$ a permutation group acting on the elements of $V$. The subset $B \subseteq V$ is a {\bf $G$-block} if for every $g \in G$,
either $g(B)=B$, or $g(B)\cap B=\emptyset$.
\end{defn}

In some cases, the group $G$ is clear from the context and we simply
refer to $B$ as a block.

It is a simple matter to realise that if $B$ is a $G$-block, then for any
$g \in G$, $g(B)$ will also be a $G$-block.  Also, intersections
of $G$-blocks remain $G$-blocks.

Let $G$ be a transitive permutation group, and let $B$ be a $G$-block.
Then, as noted above, $\{g(B):g\in G\}$ is a set of blocks that 
(since $G$ is transitive) partition the set $V$.  We call this set the
{\bf complete block system} of $G$ generated by the block $B$.

Notice that any singleton in $V$, and the entire set $V$, are always
$G$-blocks.

The {\bf size} of the block $B$ is the cardinality of the set $B$.
A block $B$ is {\bf nontrivial} if the size of $B$ is neither 1 nor the
cardinality of $V$.  

\begin{defn}
The transitive permutation group $G$ is said to be {\bf imprimitive} if $G$
admits nontrivial blocks.  If $G$ is transitive but not imprimitive, then $G$ is said to
be {\bf primitive}.
\end{defn}

The notion of transitivity for permutation groups can be generalised.
\begin{defn}
The permutation group $G$ acting on the set $V$ is {\bf $k$-transitive} if
given any two $k$-tuples $(v_1,\ldots,v_k)$ and $(u_1,\ldots, u_k)$ with
$v_1,\ldots,v_k,u_1,\ldots, u_k \in V$, there exists some $g \in G$ such
that $g(v_i)=u_i$ for $1 \le i \le k$.

In particular, we often say that a 2-transitive group is doubly transitive.
\end{defn}
\begin{defn}
The abstract group $G$ is a {\bf Burnside group} if every primitive
permutation group containing the regular representation of $G$ as
a transitive subgroup is doubly transitive.
\end{defn}

Burnside gave the first example of such a group, hence the name.  This
is extremely useful in the theory of circulant graphs, due to the
following theorem.

\begin{thm}[Theorem 25.3, \cite{W}]
Every cyclic group of composite order is a Burnside group.  
\end{thm}

In particular, if $n$ is composite, the automorphism group of a circulant graph of order $n$
is either doubly transitive or imprimitive.
\begin{cor}
For any circulant graph $X=X(n;S)$, one of the following holds:
\begin{enumerate}
\item Aut($X$)$=\sym_n$;
\item Aut($X$) is imprimitive; or
\item $n$ is prime.
\end{enumerate}
\end{cor}

This can be used repeatedly to show that the minimal blocks of a circulant
graph must either have prime size, or the induced subgraph on these vertices is complete or empty.

We can also use this to determine additional information about how Aut($X$)
acts upon any blocks - once again, this action must be imprimitive or doubly transitive unless there are a prime number of blocks.

\section{Related problems}

In this section, we discuss a number of problems that are closely related to finding the automorphism group of a circulant graph.

\subsection{The Cayley Isomorphism problem for circulants}

One of the most-studied problems related to finding the automorphism group of circulant graphs, is the Cayley Isomorphism, or CI, problem. 

\begin{defn}
A circulant graph $X=X(n;S)$ is said to have the {\bf Cayley Isomorphism (CI) property} if whenever $Y=Y(n;S')$ is isomorphic to $X$, there is some $a \in \Z_n^*$
for which $aS=S'$.

The cyclic group of order $n$ is said to have the Cayley Isomorphism (CI) property if every circulant graph $X(n;S)$ has the CI property.
\end{defn}

The CI problem, of course, is to determine which graphs (or which groups) have the CI property.  When applied to digraphs, the CI property is referred to as the DCI property (for ``directed Cayley Isomorphism" property).

\begin{thm}\cite{APa, B}
A circulant graph $X$ on $n$ vertices has the CI property if and only if any
two $n$-cycles in Aut($X$) are conjugate in Aut($X$).
\end{thm}
So, knowing that a graph has the CI-property gives us significant information about its automorphism group.

We give a brief history of the major results on this problem.

In 1967, \`{A}d\`{a}m conjectured \cite{Adam} that all cyclic groups are CI-groups.
Elspas and Turner proved in 1970 that this was not the case \cite{ET}, and that, in fact, $\Z_{p^2}$ is 
not CI for $p \ge 5$.  There were also many positive results on this conjecture, however.  First, in 1967, Turner proved that $\Z_p$ is CI \cite{Turner}.  In a computer search, McKay \cite{McKay} found that $\Z_n$ is CI, when $n\le 37$ and $n \not= 16, 24, 25, 27, 36$.  In 1977, Babai \cite{B} proved that $\Z_{2p}$ is CI, which was generalised in 1979 by Alspach and Parsons \cite{APa}, who proved that $\Z_{pq}$ is CI.  In 1983, Godsil \cite{Godsil} proved that $\Z_{4p}$ is CI.  Muzychuk completed the work in 1997 \cite{Muz1, Muz2} by showing that  
$\Z_n$ is DCI if and only if $n \in \{k,2k,4k\}$ where $k$
is odd and squarefree; $\Z_n$ is CI if and only if $n \in \{8,9,18\}$ or $\Z_n$
is DCI.

This problem has also been studied for particular families of graphs; 
Huang and Meng proved in 1996 \cite{Huang}, that $X(n;S)$ has the CI property if $S$ is a minimal generating set for $\Z_n$.
In 1977,
Toida conjectured that if $S\subseteq\Z_n^*$ then $X(n;S)$ is CI \cite{Toida}.
This conjecture was proven by Klin, P\"{o}schel and Muzychuk \cite{Muzy2}, and independently by Dobson and Morris \cite{ourToida}.

The special case where the connection set, $S$, is small, has also been studied.
In 1977, Toida \cite{Toida} proved that if $|S|\le 3$, $X(n;S)$ is CI.
In 1988, Sun \cite{Sun} was the first to prove that if $|S|=4$, $X(n;S)$ is CI, although others later also proved this result. In 1995, Li proved that if $|S|=5$, $X(n;S)$ is CI \cite{Li1}.

This is by no means intended as a full history of work that has been done on the CI problem;
for further information, the reader is referred to Li's survey of the problem \cite{Lisurvey}.

\subsection{Edge-transitivity, arc-transitivity, 2-arc-transitivity}

\begin{defn}
A $k$-arc is a list $v_1, \ldots, v_k$ of vertices for which
any two sequential vertices are adjacent, and any 3 sequential vertices are distinct.
\end{defn}

This definition allows us to consider graphs whose automorphism groups are transitive on the set of $k$-arcs of the graph, for various values of $k$; for our purposes, we will only consider $k=1$ and $k=2$.

The following theorem was proven by Chao \cite{Chao} in 1971; the proof was simplified by Berggren \cite{Berggren} in 1972.
\begin{thm}
A circulant graph on $p$ vertices is edge-transitive iff $S=\emptyset$ or
$S$ is a coset of an even order subgroup $H \le \Z_p^*$.
\end{thm}
Notice that in the case of circulant graphs, the reflection is an automorphism, so arc-transitivity is equivalent to edge-transitivity. 

This result, and those that follow, provide significant information about how the automorphism group of a graph can act on the edge set of the graph. They are therefore related to finding the automorphism group, although they do not directly provide much information about the automorphism group of an arbitrary circulant graph.

The above theorem was extended to classifications of arc-transitive graphs on $pq$ vertices, for any distinct primes $p$ and $q$, in papers by Cheng and Oxley \cite{ChOx}, who classified the graphs on $2p$ vertices; Wang and Xu \cite{WangXu}, who classified the graphs on $3p$ vertices; and Praeger, Wang and Xu \cite{PrWaXu}, who completed the classification.

In 2001, the following classification was obtained of arc-transitive circulants on a square-free number of vertices, by Li, Maru\v{s}i\v{c} and Morris \cite{arc-trans}.
\begin{thm}
If $X$ is an arc-transitive circulant graph of square-free order $n$, then one of the following holds:
\begin{enumerate}
\item $X=K_n$;
\item Aut($X$) contains a cyclic regular normal subgroup; or
\item $X=Y\wr \bar{K_b}$, or $X=Y\wr \bar{K_b}-bY$, where $n=mb$, and $Y$ is
an arc-transitive circulant of order $m$.
\end{enumerate}
\end{thm}

Another approach has been taken to classifying edge-transitive circulant graphs, using the more stringent condition that the complement must also be edge-transitive.  The following result was proven by Zhang in 1996 \cite{Zhang}.

\begin{thm}
If $G$ and $\overline{G}$ are both edge-transitive circulants, then $G$ is
one of: $mK_n$, $\overline{mK_n}$, or a self-complementary Paley graph
($n$ is prime, 1 (mod 4) and $S=\{a^2: a \in \Z_n^*\}$).
\end{thm}

The property of 2-arc-transitivity is stronger than arc-transitivity, but accordingly tells us more about the action of the automorphism group.  Circulant graphs that are 2-arc-transitive have been fully classified by Alspach, Conder, Maru\v{s}i\v{c} and Xu \cite{ACMX}, as follows.
\begin{thm}
A connected 2-arc-transitive circulant graph is one of:
\begin{enumerate}
\item $K_n$ (exactly 2-arc-transitive);
\item $K_{n/2,n/2}$ (exactly 3-arc-transitive);
\item $K_{n/2,n/2}$ minus a 1-factor, $n\ge 10$, $n/2$ odd (exactly 2-arc-transitive;
\item $C_n$ ($k$-arc-transitive for every $k \ge 0$)
\end{enumerate}
\end{thm}

\subsection{Other regular subgroups in Aut($X$)}

The automorphism group of any circulant graph will have a cyclic subgroup that acts regularly on the vertices of the graph. (A permutation group on a set $V$ is said to act regularly if for any pair of points in $V$, there is exactly one permutation in the group that maps one to the other.)  Sometimes, the automorphism group of a graph may have multiple, nonisomorphic, regular subgroups.  This is of interest from a different perspective, because it means that the graph in question can be represented as a Cayley graph on some noncyclic group, besides being a circulant graph.  From our perspective, knowing the regular subgroups of the automorphism group may be useful in determining the automorphism group.

There are only a few results of note on this topic.  The first result was proven by Joseph in the special case $n=p^2$ \cite{Joseph}, and extended by Morris \cite{different} to all prime powers.
\begin{thm}
If $n=p^e$, $X$ is a circulant graph on $n$ vertices, and Aut($X$) 
contains a regular subgroup that is not cyclic, then $X$ is isomorphic to a wreath product of smaller circulant graphs.
\end{thm}

Recently, Maru\v{s}i\v{c} and Morris \cite{Dragan2000} proved the following results.
\begin{thm}
Let $X=X(n;S)$ be a circulant graph, and $\Z_n^*(S)$ be the subgroup of $\Z_n^*$ that fixes $S$ set-wise.  Then if $\gcd(n,|\Z_n^*(S)|)>1$,
the automorphism group of $X$ has a noncyclic regular subgroup.  

In fact, if $p$ is any prime divisor of $\gcd(n,|\Z_n^*(S)|)$, then the automorphism group of $X$ contains a regular subgroup that is isomorphic to $\Z_p\ltimes \Z_{n/p}$.
\end{thm}
They also show that the converse is not true in general; that is, there may be noncyclic regular subgroups in the automorphism group even if this greatest common divisor is 1.

In the same paper, they prove the following result.  The condition that we must have a normal circulant means that the regular cyclic group must be normal in the automorphism group of the graph.
\begin{thm}
Let $X$ be a normal circulant graph of order $n$, $n$ not divisible by 4.  Then if the automorphism group of $X$ has a noncyclic regular subgroup, that group must be metacyclic, generated by two cyclic subgroups whose orders are relatively prime.
\end{thm}

\section{Concluding remarks}

Efficient algorithms have been found for determining the automorphism group of a circulant graph, when the number of vertices is any product of distinct primes. In the case where the number of vertices is a prime power, structural theorems about the automorphism groups exist (in unpublished form), but no algorithms have been constructed.  No work has been done on combining these into results that may hold when the number of vertices is divisible by at least two distinct primes, but is not square-free.  There is, therefore, a great deal of room for more results on this problem.

Unfortunately, in the results that we have seen, the complexity (of the proofs) seems to be growing.  The result on circulant graphs on a square-free number of vertices is the culmination of two densely-packed, long papers; the first provides a structural theorem, and the second the algorithm.  The result that deals with circulant graphs whose number of vertices is a prime power seems to be unpublished to this point, largely because of its length and complexity.  It may be that new techniques will need to be developed before the problem can be completed in its full generality.

Of course, circulant graphs are just one very small step towards answering the general question with which we began, of finding the automorphism group of any graph.  In full generality, this problem seems unmanageable, but it may be that for some other classes of graphs, the solution turns out to be feasible, or even easy.

\end{document}